\documentclass[11pt]{article}
\usepackage{amsfonts}
\usepackage{euscript,amstext}
\usepackage{mathrsfs}
\usepackage{indentfirst}
\usepackage{enumerate}
\usepackage{amssymb}
\usepackage{leftidx}
\usepackage{amsmath}
\usepackage{hyperref}
\usepackage{amscd}
\numberwithin{equation}{section}
\newtheorem{thm}{Theorem}[section]
\numberwithin{thm}{section}

\newtheorem{prop}[thm]{Proposition}
\newtheorem{rem}[thm]{Remark}

\newtheorem{cor}[thm]{Corollary}

\DeclareMathOperator{\tr}{Tr}

\DeclareMathOperator{\dist}{d}

\DeclareMathOperator{\im}{im}
\DeclareMathOperator{\dom}{dom}

\DeclareMathOperator{\res}{Res}

\title{Comparison between two analytic torsions
\\
on orbifolds}
\author{Xianzhe DAI\footnote{Department of Mathematics, University of California, Santa Barbara,
CA 93106, USA (dai@math.ucsb.edu)}\ \ and\  \
Jianqing YU\footnote{School of Mathematical Sciences,
University of Science and Technology of China,
96 Jinzhai Road, Hefei, Anhui 230026,
P. R. China. (jianqing@ustc.edu.cn)}
}

\begin{document}

\maketitle

\begin{abstract}
In this paper, we establish an equality between the analytic torsion
introduced by Dar\cite{MR876230} and the orbifold analytic torsion defined by Ma \cite{MR2140438}
on a compact manifold with
isolated conical singularities which in addition has an orbifold structure.
We assume the orbifold flat vector bundle is an honest vector bundle, although the metric on the flat bundle may not be flat.
\end{abstract}

\section{Introduction}
Let $X$ be a closed Riemannian manifold and $F$ be a flat real vector bundle over $X$ carrying the flat connection $\nabla^F$. As a
geometric invariant, the analytic torsion of $X$, was first introduced by Ray and Singer \cite{MR0295381,MR0383463}
in searching for an analytic interpretation of the combinatorially defined Reidemeister torsion. The latter is
the first topological invariant which is not a homotopy invariant, introduced by
Reidemeister \cite{Reidemeister1935} and Franz \cite{franz1935torsion}.

The Ray-Singer analytic torsion is a weighted linear combination of the
determinants of the Laplacian acting on the space of differential forms with values in $F$, and depends on
the the metrics on $F$ and on the Riemannian manifold. We explain this in greater detail.
Let $g^{TX}$ be a Riemannian metric on $X$ and $g^F$ be an Euclidean metric on $F$. We denote
by $\Omega^i(X,F)$ the space of smooth $i$-forms on $X$ with values in $F$, and
set $\Omega(X,F)=\oplus_{i}\Omega^i(X,F)$. The flat connection $\nabla^F$ induces naturally
a differential $d^F$ on $\Omega(X,F)$. Let $d^{F*}$ be the (formal) adjoint operator of $d^F$
associated with $g^{TX}$ and $g^F$. The Laplacian operator $\Delta:=d^Fd^{F*}+d^{F*}d^F$
acts on $\Omega(X,F)$ and preserves its $\mathbb{Z}$-grading. Let $P$ be the orthogonal projection
operator from $\Omega(X,F)$ to $\ker \Delta$.
Let $e^{-t\Delta}$ be the heat semi-group operator of $\Delta$.

Let $\Gamma$ be the gamma function. The torsion zeta function is defined as follows.
For $s\in \mathbb{C}$, $\Re s>\tfrac{\dim X}{2}$,
\begin{equation}\label{a8}
\zeta(s)=\frac{-1}{\Gamma(s)}\sum_{i=0}^{\dim X}(-1)^ii\int^\infty_0
t^{s-1}\tr\big|_{\Omega^i(X,F)}\big[e^{-t\Delta}(1-P)\big]dt.
\end{equation}
By the standard elliptic theory and the asymptotic expansion of heat kernel, $\zeta(s)$ extends to a meromorphic
function of $s\in \mathbb{C}$ which is holomorphic at $s=0$.
The Ray-Singer analytic torsion of $X$ with coefficients in $F$ is defined as
\begin{equation}\label{a7}
T(X,g^{TX},g^F)=\exp\left\{ \frac{1}{2}\frac{\partial \zeta}{\partial s}(0)\right\}.
\end{equation}

The celebrated Cheeger-M\"uller theorem \cite{MR528965,MR498252} establishes
the Ray-Singer conjecture. Namely,
if $g^F$ is flat, i.e. $(F,\nabla^F)$ is induced by
an orthogonal representation of the fundamental group of $X$,
then the analytic torsion coincides with the classical Reidemeister
torsion. There are various generalizations of this result.
M\"uller \cite{MR1189689} extended his result to the case of unimodular representation,
where only the metric induced on $\det F$ is required to be flat.
Bismut and Zhang \cite{MR1185803} reformulated the above Cheeger-M\"uller theorem as an equality
between the Reidemeister and Ray-Singer metrics defined on the determinant of cohomology, and proved an extension of
it to arbitrary flat vector bundles with arbitrary Euclidean metrics.

Based on the work of Goresky, MacPherson \cite{MR572580,MR696691} and Cheeger \cite{MR530173,MR573430,MR730920},
Dar \cite{MR876230} described a possible extension of the Cheeger-M\"uller theorem to singular spaces.
Assume now that $X$ is a Riemannian  manifold with isolated conical singularities $\Sigma$
and $F$ be a flat bundle over $X\backslash\Sigma$.
Let $\Omega^i_0(X\backslash\Sigma,F)$ be the space of smooth $i$-forms with values in $F$, compactly supported on $X\backslash\Sigma$, and
set $\Omega_0(X\backslash\Sigma,F)=\oplus_{i}\Omega^i_0(X\backslash\Sigma,F)$.
Choose an ideal boundary condition in the sense of Cheeger \cite{MR530173,MR573430}, which
corresponds to a closed extension of the
de Rham complex $(\Omega_0(X\backslash\Sigma,F),d^F)$ into a Hilbert complex.
Let $\Delta_{\rm c}$ be the Laplacian associated to the Hilbert complex,
and $e^{-t\Delta_{\rm c}}$ be the heat semi-group operator of $\Delta_{\rm c}$.

Although the asymptotic expansion of the trace of the heat kernel of $e^{-t\Delta_{\rm c}}$
can contain logarithmic terms \cite{MR730920}, Dar \cite{MR876230} made a crucial observation that
after taking the weighted linear combination as in \eqref{a8},  the contribution from
the logarithmic term drops out.
Thus, one can define the analytic torsion $T_{\rm c}(X,g^{TX},g^F)$ as in \eqref{a7}.

On the other hand, Ma \cite{MR2140438} extended the Quillen metric to compact complex orbifolds. In this setting, he
established the anomaly formula and
calculated the behavior of the Quillen metric by orbifold immersions,
which generalized the corresponding results in \cite{MR931666} and \cite{MR1188532}.

We assume in addition that $X$ carries an orbifold structure and the relevant geometric data are all in the
orbifold category. Following Ma \cite{MR2140438}, one can define the orbifold torsion $T_{\rm o}(X,g^{TX},g^F)$
in this real setting.

The main result of this paper is an equality between these two analytic torsions on a special type of singular space,
which are defined from the rather distinct perspectives to the singularities.
\begin{thm}\label{main}
Let $X$ be a compact manifold with
isolated conical singularities which has an orbifold structure
and $F$ be a flat real orbifold vector bundle over $X$ with trivial isotropic action on its fibers.
For any (orbifold) Riemannian metric $g^{TX}$ on $X$ and
(orbifold) Euclidean metric $g^F$ on $F$, one has
\begin{equation}
T_{\rm c}(X,g^{TX},g^F)= T_{\rm o}(X,g^{TX},g^{F}).
\end{equation}
\end{thm}

The proof of Theorem \ref{main} consists of the comparison of the two kinds of cohomology and of the two heat kernels.

The isomorphism between the two kinds of cohomology can be established by a straightforward calculation using the Mayer Vietoris sequence
and Cheeger's result \cite{MR530173} on $L^2$ cohomology.

The comparison of the two heat kernels is inspired by Cheeger's approach \cite{MR528965} to the Ray-Singer
conjecture. Our basic observation is that Cheeger's constructions of the Green kernels on the annulus $A^{m+1}_{u,1}$ and the corresponding
estimates in \cite[Section 6]{MR528965} can be used with little modification to establish a Sobolev inequality
on $C_{(0,1]}(S^m/G)$, from which we get the estimate for the heat kernel constructed from the point of view of conical singularities.
Combining this with the estimate for orbifold heat kernel \cite{MR2140438} and applying the Duhamel principle \cite[(3.9)]{MR528965},
we obtain an equality between these two heat kernels outside the singularities.

The rest of this paper is organized as follows. In Section \ref{s2}, we
recall the basic definitions of the two analytic torsions
defined by Dar \cite{MR876230} and Ma \cite{MR2140438}, respectively.
Section \ref{s3} is devoted to a proof of Theorem \ref{main}.

\vspace{3mm}\textbf{Acknowledgements}\ \ The first author is supported by NSF-DMS.
The work was carried out while the second author
was visiting the University of California, Santa Barbara (UCSB).
He would like to thank the hospitality of the Department of Mathematics in UCSB and
the financial support from the program of China Scholarships Council.

\section{Two analytic torsions}\label{s2}
In this section, we first recall the relevant facts about
$L^2$ cohomology and Dar's definition of analytic torsion \cite{MR1174159} on
a manifold with isolated conical singularities. Then we introduce a real analogue of
Ma's analytic torsion \cite{MR2140438} which is originally defined on compact complex orbifolds.

\subsection{The analytic torsion on a manifold with conical singularities}

Let $N$ be a closed manifold carrying a Riemannian metric $g^{TN}$. The model cone $C(N)$ over $N$
is the space $(0,+\infty)\times N$ carrying the conical metric $dr^2+r^2 g^{TN}$, where $r$ denotes the
radial coordinate.
For $u>0$, set
\begin{equation}
\begin{split}
C_{(0,u)}(N)&=\{(r,y)\in C(N)\,|\,0<r<u \},
\\
C_{(0,u]}(N)&=\{(r,y)\in C(N)\,|\,0<r\leq u \}.
\end{split}
\end{equation}

Let $X$ be a $m+1$ dimensional compact Riemannian manifold with isolated conical singularities $\Sigma$.
By this we mean that $X\backslash\Sigma$ is a smooth manifold of dimension $m+1$ with a Riemannian metric $g^{TX}$,
and for $p\in \Sigma$, there exists an open neighborhood $\mathscr{U}_u(p)$ ($u>0$) in $X$ and a
closed Riemannian manifold $(N_p,g^{TN_p})$
such that $(\mathscr{U}_u(p)\backslash\{p\},g^{TX}|_{\,\mathscr{U}_u(p)\backslash\{p\}})$ is isometric to
the cone $(C_{(0,u)}(N_p), dr^2+r^2 g^{TN_p})$.
For $u>0$ small enough, set
\begin{equation}
C^*_{(0,u)}(N_p)=C_{(0,u)}(N_p)\cup\{p\},\quad C^*_{(0,u]}(N_p)=C_{(0,u]}(N_p)\cup\{p\}.
\end{equation}
Then $X_u:=X\backslash\big(\sqcup_{p\in\Sigma}C^*_{(0,u)}(N_p)\big)$ is a compact
manifold with boundary,
\begin{equation}\label{t3}
X=X_u\cup(\mathop{\sqcup}_{p\in \Sigma}C^*_{(0,u]}(N_p)\,),
\end{equation}
where the union is along the boundary $\sqcup_{p\in \Sigma}\{u\}\times N_p$.

Let $F$ be a flat real vector bundle over $X\backslash\Sigma$ carrying the flat connection $\nabla^F$.
Let $\Omega^i_0(X\backslash\Sigma,F)$ be the space of smooth $i$-forms with values in $F$, compactly supported on $X\backslash\Sigma$, and
set $\Omega_0(X\backslash\Sigma,F)=\oplus_{i}\Omega^i_0(X\backslash\Sigma,F)$.
Let $L^{2}(\Omega^i(X\backslash\Sigma,F))$ be the Hilbert space which consists of
square integrable $i$-forms on $X\backslash\Sigma$ with values in $F$.

We recall the main features of the $L^2$-cohomology of X with coefficients in $F$.
Let $(\Omega_0(X\backslash\Sigma,F),d^F)$ be the de Rham complex, where $d^F$ is the differential induced from
the flat connection $\nabla^F$ in a natural way.
Recall that
an ideal boundary condition (cf. \cite[Section 2]{MR573430}, \cite[p. 105]{MR1174159}) for the complex $(\Omega_0(X\backslash\Sigma,F),d^F)$
is a choice of closed extensions $d^F_{{\rm c},i}$ of $d^F_i$ in $L^{2}(\Omega^i(X\backslash\Sigma,F))$ for $0\leq i\leq m+1$,
such that
\begin{equation}
d^F_{{\rm c},i}(\dom(d^F_{{\rm c},i}))\subset\dom(d^F_{{\rm c},i+1}),\quad d^F_{{\rm c},i+1}\circ d^F_{{\rm c},i}=0.
\end{equation}
We then get a Hilbert complex in the sense of Br\"uning-Lesch \cite[p. 90]{MR1174159},
\begin{equation}\label{z1}
0\rightarrow\dom(d^F_{{\rm c},0})\xrightarrow[]{d^F_{{\rm c},0}}\cdots\xrightarrow[]{d^F_{{\rm c},m}}\dom(d^F_{{\rm c}.m+1})\rightarrow 0.
\end{equation}
Let $\delta^F$ denote the formal adjoint of $d^F$, then the minimal and maximal extensions of $d^F$,
\begin{equation*}
d^F_{\min}:=\text{closure of }d^F, \quad d^F_{\max}:=\text{adjoint of the closure of }\delta^F,
\end{equation*}
are examples of the ideal boundary conditions. A prior there may be several distinct ideal boundary conditions.

We assume that 
\begin{equation}
H^{\frac{m}{2}}(N_p,F|_{N_p})=0,\ p\in \Sigma.
\end{equation}
As shown in \cite{MR573430}, in this case the ideal boundary condition is unique, i.e.,
\begin{equation}
d^F_{\min,i}=d^F_{\max,i}\,,  \quad\text{for}\ 0\leq i\leq m+1.
\end{equation}
We denote by $(\mathcal{C},d_{\rm m})$ the unique extension of $(\Omega_0(X\backslash\Sigma,F),d^F)$
into a Hilbert complex as in \eqref{z1}. The $L^2$-cohomology of X with coefficients in $F$ is
the cohomology of the complex $(\mathcal{C},d_{\rm m})$,
\begin{equation}
H^i_{(2)}(X,F):= \ker d_{{\rm m},i}/\im d_{{\rm m},i-1}\,, \quad 0\leq i\leq m+1.
\end{equation}

Let $\Delta_{\rm c}$ be the Laplacian associated to the Hilbert complex $(\mathcal{C},d_{\rm m})$,
\begin{equation}
\Delta_{{\rm c}}:=d^F_{\min}\delta^F_{\min}+\delta^F_{\min}d^F_{\min}\,.
\end{equation}
We denote by $\Delta_{{\rm c},i}$ the restriction of $\Delta_{\rm c}$ to $\mathcal{C}^i$ for $0\leq i\leq m+1$.
By the $L^2$-Hodge theorem \cite[Section 1]{MR573430}, one knows
all the $L^2$-cohomology groups $H^i_{(2)}(X,F)$ are finite dimensional
and the complex $(\mathcal{C},d_{\rm m})$ is Fredholm in the sense of
Br\"uning-Lesch \cite[p. 90]{MR1174159}. Moreover, the canonical maps
\begin{equation}\label{c1}
\ker\Delta_{{\rm c},i}\longrightarrow H^i_{(2)}(X,F), \quad 0\leq i\leq m+1.
\end{equation}
are isomorphisms (see also \cite[Corollary 2.5]{MR1174159}).

Let $e^{-t\Delta_{{\rm c},i}}$ be the heat semi-group operator of $\Delta_{{\rm c},i}$, and
$K_{{\rm c},i}(t,\cdot,\cdot)$ be the heat kernel of $e^{-t\Delta_{{\rm c},i}}$.
Let $P_{{\rm c},i}$ be the orthogonal projection
operator from $L^{2}(\Omega^i(X\backslash\Sigma,F))$ to $\ker \Delta_{{\rm c},i}$.
Set $P^\perp_{{\rm c},i}=1-P_{{\rm c},i}$.

One defines the $i$-th torsion zeta function as follows.
\begin{equation}
\zeta_{{\rm c},i}(s):=\frac{-1}{\Gamma(s)}\int^\infty_0t^{s-1}\tr\big[e^{-t \Delta_{{\rm c},i}}P^\perp_{{\rm c},i}\big]dt.
\end{equation}
The asymptotic expansion of the trace of the heat kernel $K_{{\rm c},i}(t,\cdot,\cdot)$ yields a meromorphic extension of
$\zeta_{{\rm c},i}(s)$ to the whole complex plane, and determines its behavior near $s=0$.
In particular, it need not be regular at $s=0$ because of the appearance of a logarithmic term \cite[Theorem 2.1]{MR730920}.
The crucial observation made by Dar \cite[Theorem 4.4]{MR876230} is that
\begin{equation}
\mathop{\res}_{s=0}\Big(\sum_{i=0}^{m+1}(-1)^ii\cdot\zeta_{{\rm c},i}(s)\Big)=0.
\end{equation}

Thus, the full torsion zeta function $\zeta_{\rm c}(s):=\sum_{i=0}^{m+1}(-1)^ii\cdot\zeta_{{\rm c},i}(s)$ is
indeed holomorphic at $s=0$. The analytic torsion $T_{\rm c}(X,g^{TX},g^F)$ of $X$ with coefficients in $F$ is defined as in \eqref{a7}
(cf. \cite[p. 215]{MR876230}),
\begin{equation}\label{t4}
T_{\rm c}(X,g^{TX},g^F)=\exp\left\{\frac{1}{2}\frac{\partial \zeta_{\rm c}}{\partial s}(0)\right\}.
\end{equation}

\subsection{The analytic torsion on an orbifold}

We refer to \cite{MR0474432} for relevant definitions of orbifolds
and to \cite{MR2140438} for notations used here. In \cite{MR0474432} orbifolds were called V-manifolds.

Let $(X,\mathcal{U})$ be a compact orbifold endowed with a Riemannian metric $g^{TX}$.
Let $F$ be a flat real orbifold vector bundle over $X$ equipped with the flat connection $\nabla^F$
and an Euclidean metric $g^F$.

We denote by $\Omega^i(X,F)$ the space of smooth sections of $\Lambda^i(T^*X)\otimes F$ over $X$, and
set $\Omega(X,F)=\oplus_{i}\Omega^i(X,F)$.
The flat connection $\nabla^F$ induces naturally
a differential $d^F$ on $\Omega(X,F)$.
Let $H(X,F)=\oplus_{i=0}^{\dim X}H^i(X,F)$ be the singular cohomology group of $X$ with
coefficients in $F$. The de Rham theorem for orbifolds \cite[p. 78]{MR0474432} gives
us a canonical isomorphism,
\begin{equation}\label{b7}
H^i(\Omega(X,F),d^F)\simeq H^i(X,F),\ \text{for}\ 0\leq i\leq \dim X.
\end{equation}

Let $\langle\cdot,\cdot\rangle_{\Lambda(T^*X)\otimes F}$ be the
metric on $\Lambda(T^*X)\otimes F$ induced from $g^{TX}$, $g^F$,
and $dv_X$ be the Riemannian volume form on $X$ associated to $g^{TX}$.
As in \cite[(2.8)]{MR2140438}, one defines the $L^2$-scalar product
$\langle\cdot,\cdot\rangle$ on $\Omega(X,F)$, for $s,s'\in\Omega(X,F)$,
\begin{equation}\label{a9}
\langle s,s'\rangle:=\int_X\langle s,s'\rangle_{\Lambda(T^*X)\otimes F}(x)dv_X(x).
\end{equation}
Let $\delta_{\rm o}^{F}$ be the (formal) adjoint operator of $d^F$
with respect to \eqref{a9}.
Set
\begin{equation}
\Delta_{\rm o}:=d^F\delta_{\rm o}^{F}+\delta_{\rm o}^{F}d^F.
\end{equation}
Then $\Delta_{\rm o}$ is a second order differential operator, which acts
on $\Omega(X,F)$ and preserves its $\mathbb{Z}$-grading, with $\sigma_2(\Delta_{\rm o})=|\xi|^2$ ($\xi\in T^*X$).

Using the same proof as in \cite[Proposition 2.2]{MR2140438},
one deduces the Hodge decomposition.
\begin{prop}
There is an $L^2$-orthogonal direct sum decomposition,
\begin{equation}\label{b8}
\Omega^i(X,F)=\ker \Delta_{{\rm o},i}\oplus \im d_{i-1}^F \oplus \im \delta_{{\rm o},i+1}^{F},\ \text{for}\ 0\leq i\leq \dim X.
\end{equation}
\end{prop}

From \eqref{b7} and \eqref{b8}, one has the canonical identification
\begin{equation}\label{c2}
\ker \Delta_{{\rm o},i}\simeq H^i(X,F),\ \text{for}\ 0\leq i\leq \dim X.
\end{equation}

Let $K_{\rm o}(t,\cdot,\cdot)$ be the heat kernel of the heat semigroup operator $e^{-t\Delta_{\rm o}}$ with respect to
$dv_X$.
\begin{prop}{\rm (cf. \cite[Proposition 2.1]{MR2140438})}\label{es2}
For each $U\in\mathcal{U}$, there exists a smooth section
$\Phi_j\in \Gamma(\widetilde{U}\times \widetilde{U},{\rm pr}_1^*\widetilde{F}\otimes{\rm pr}_2^*\widetilde{F})$ such that
for every $k>\dim X$, $x,y\in U$, as $t\rightarrow 0$,
\begin{equation}
\begin{split}
K_{\rm o}(t,x,y)
&=\frac{(4\pi t)^{-\frac{\dim X}{2}}}{|K^F_U|}\sum_{g\in G_{U}^F}\sum_{j=0}^k
e^{-\frac{\widetilde{\dist}^2(g\widetilde{x},\widetilde{y})}{4t}} g^{-1}\,\Phi_j(g\widetilde{x},\widetilde{y})\,t^j
\\
&\hspace{130pt}+O(t^{k-\frac{\dim X}{2}}).
\end{split}
\end{equation}
On $\{(x,y)\in X\times X\,|\,\dist^X(x,y)>c>0\}$, we have
\begin{equation}\label{b9}
\partial^\alpha_x\partial^\beta_yK_{\rm o}(t,x,y)=O(e^{-\frac{c^2}{4t}}),\ \text{as}\ t\rightarrow 0.
\end{equation}
\end{prop}

From \eqref{b9}, one has the following estimate for the pointwise norm of ${K}_{\rm o}(t,\cdot,\cdot)$.
\begin{cor}
Given $T>0$, $\varepsilon>0$ and $n\in \mathbb{N}$, there exists a constant $K(T,\varepsilon,n)>0$ such that
for $x, y\in X$ with $\dist^X(x,y)>\varepsilon$, and $0<t\leq T$,
\begin{equation}
\|{K}_{\rm o}(t,x,y)\|\leq K(T,\varepsilon,n)\,t^n.
\end{equation}
The same estimates hold for $d_x^F{K}_{\rm o}(t,x,y)$ and $\delta_{{\rm o},y}^F {K}_{\rm o}(t,x,y)$.
\end{cor}

Let $P_{\rm o}$ be the orthogonal projection operator from $\Omega(X,F)$ on $\ker \Delta_{\rm o}$ with
respect to the $L^2$ scalar product. Set $P^\perp_{\rm o}=1-P_{\rm o}$.

For $s\in \mathbb{C}$, $\Re s>\tfrac{\dim X}{2}$, set
\begin{equation}
\zeta_{\rm o}(s)=\frac{-1}{\Gamma(s)}\sum_{i=0}^{\dim X}(-1)^ii\int^\infty_0
t^{s-1}\tr\big|_{\Omega^i(X,F)}\big[e^{-t\Delta_{\rm o}}P^\perp_{\rm o}\big]dt.
\end{equation}
Using Proposition \ref{es2}, $\zeta_{\rm o}(s)$ extends to a meromorphic
function of $s\in \mathbb{C}$ which is holomorphic at $s=0$.
The orbifold analytic torsion of $(X,\mathcal{U})$ with coefficients in $F$ is defined as in \eqref{a7},
\begin{equation}
T_{\rm o}(X,g^{TX},g^F)=\exp\left\{ \frac{1}{2}\frac{\partial \zeta_{\rm o}}{\partial s}(0)\right\}.
\end{equation}

\section{The equality between two torsions}\label{s3}

In this section, we first prove a Sobolev inequality on the bounded cone $C_{(0,1]}(S^m/G)$. Then for an even
dimensional manifold with isolated conical singularities which carries an orbifold structure,
we establish an isomorphism between the $L^2$ cohomology and the singular cohomology and
an equality between $K_{\rm c}(t,\cdot,\cdot)$ and $K_{\rm o}(t,\cdot,\cdot)$.
As a corollary, we deduce our main result.

\subsection{The Sobolev inequality on a bounded cone}
Let $N=S^m/G$, the quotient space induced from
a free action of a finite group $G$ on $S^{m}$, which is a
closed orientable manifold.
Fix a Riemannian metric $g^{TN}$ on $N$.
In this subsection, we establish the Sobolev inequality on the bounded cone
$C_{(0,1]}(N)=\{(r,y)\in C(N)\,|\,0<r\leq 1\}$.

Let $F$ be a flat real vector bundle over $C(N)$ with the flat connection $\nabla^F$ and an Euclidean metric $g^F$.
By parallel transport along the radial geodesic with respect to $\nabla^F$,
one can identify $F|_{\{r\}\times N}$ with $F|_{\{1\}\times N}$,
and $g^{F|_{\{r\}\times N}}$ with $g^{F|_{\{1\}\times N}}$.
Let $\pi: C(N)\longrightarrow \{1\}\times N$ given by $(r,y)\longmapsto (1,y).$
Then we have the following identification
\begin{equation}\label{a1}
(F,\nabla^F,g^F)=\pi^*\big(F|_{\{1\}\times N},\nabla^{F|_{\{1\}\times N}},g^{F|_{\{1\}\times N}}\big),
\end{equation}
and view $(F|_{\{1\}\times N},\nabla^{F|_{\{1\}\times N}},g^{F|_{\{1\}\times N}})$ as a flat vector bundle over
$N$ by identifying $\{1\}\times N$ with $N$.

Let $F^*$ be the dual bundle of $F$ carrying the (dual) flat connection $\nabla^{F^*}$ and the (dual) metric $g^{F^*}$.

Let $\Omega^i(C(N),F)$ and $\Omega^i(C(N),F^*)$ denote the spaces of smooth $i$-forms on $C(N)$ with values in $F$
and $F^*$ respectively. Set
\begin{equation*}
\Omega(C(N),F)=\oplus_{i=0}^{m+1}\Omega^i(C(N),F), \ \Omega(C(N),F^*)=\oplus_{i=0}^{m+1}\Omega^i(C(N),F^*).
\end{equation*}

Let $d^F$ and $d^{F^*}$ be be the exterior differentials on $\Omega(C(N),F)$ induced by $\nabla^F$
and on $\Omega(C(N),F^*)$ induced by $\nabla^{F^*}$ respectively. We denote by $\delta^F$
the formal adjoint of $d^F$ with respect to the natural $L^2$ metric on $\Omega(C(N),F)$ induced by $g^{TN}$ and $g^F$. Set
\begin{equation}\label{g1}
\Delta=d^F\delta^F+\delta^Fd^F.
\end{equation}

We choose $dr\wedge r^{m} dv_{N}$ as the oriented volume form of $C(N)$,
where $dv_{N}$ is the oriented volume form of $N$. Set
\begin{equation*}
*^F: \Omega^i(C(N),F)\longrightarrow \Omega^{m+1-i}(C(N),F^*), \quad \beta\longmapsto \langle\cdot,*\beta\rangle_{g^F},
\end{equation*}
where $*$ is the usual Hodge star operator. Then one has
\begin{equation}
\delta^F|_{\Omega^i(C(N),F)}=(-1)^{i}(*^F)^{-1}d^{F^*}*^F.
\end{equation}

Let $\Omega^i(N,F)$ denote the space of smooth $i$-forms on $N$ with values in
$F|_{N}$. Set $\Omega(N,F)=\oplus_{i=0}^{m}\Omega^i(N,F)$.
As in \cite[Section 3]{MR730920}, operations on the cross section $N$
are indicated by a tilde. Let $\widetilde{d^F}$ be natural exterior differential on $\Omega(N,F)$
with formal adjoint $\widetilde{\delta^F}$. Set
\begin{equation}\label{g2}
\widetilde{\Delta}=\widetilde{d^F}\widetilde{\delta^F}+\widetilde{\delta^F}\widetilde{d^F}.
\end{equation}

From \eqref{g1}-\eqref{g2}, it is a straightforward calculation to show that for
\begin{equation*}
\beta=g(r)\phi+f(r)dr\wedge \psi\in \Omega^i(C(N),F),
\end{equation*}
where $\phi\in \Omega^i(N,F)$, and $\psi\in\Omega^{i-1}(N,F)$, one has (compare with \cite[(6.2)]{MR528965})
\begin{align}\label{a3}
\Delta\beta
&=(-g''-(m-2i)r^{-1}g'\,)\,\phi+r^{-2}g\,\widetilde{\Delta}\phi
-2r^{-3}g\,dr\wedge \widetilde{\delta^F}\phi
\notag\\
&\quad+(-f''-(m-2i+2)r^{-1}f'+(m-2i+2)r^{-2}f''\,)\,dr\wedge\psi
\notag\\
&\qquad+r^{-2}f\,dr\widetilde{\Delta}\psi
-2r^{-1}f\,\widetilde{d^F}\psi.
\end{align}

We make the assumption that $g^{TN}$ can be lifted to a $G$-equivariant Riemannian metric on $S^m$
and $F|_{N}$ can be lifted to a $G$-equivariant trivial flat bundle $\mathcal{F}$ on $S^m$ such that the $G$-action is trivial along the fiber.
Without loss of generality, we may assume ${\rm rk}\,\mathcal{F}=1$.

For $0\leq i\leq m$, we choose an orthonormal basis $\{\phi^i_j\}_{j=1}^{\infty}$
of the space $\ker\widetilde{\delta^F}\cap\Omega^{i}(N,F)$, such that
\begin{equation}
\widetilde{\Delta} \phi^i_j=\mu_j(i)\phi^i_j, \ \text{with}\ 0\leq\mu_1(i)\leq \mu_2(i)\leq\cdots.
\end{equation}

We use the same simplified notations as in \cite[(3.1)-(3.3)]{MR730920},
\begin{equation}\label{b3}
\alpha(i)=\tfrac{1+2i-m}{2},\quad \nu_j(i)=\sqrt{\mu_j(i)+\alpha^2(i)},\quad a_j^{\pm}(i)=\alpha(i)\pm\nu_j(i).
\end{equation}

We now proceed as in  Cheeger \cite[p. 289-291]{MR528965} to construct the Green operators in various cases.
\begin{enumerate}[{\rm (I)}]
\item For $\nu_j(i)>0$, set (cf. \cite[(6.6)]{MR528965})
\begin{equation}
h_{\mu_j(i)}(r_1,r_2)=\frac{1}{2\nu_j(i)}\cdot
\begin{cases}\label{b4}
r_1^{a_j^+(i)}r_2^{a_j^-(i)}, & r_1\leq r_2,
\\
r_1^{a_j^-(i)}r_2^{a_j^+(i)}, & r_2\leq r_1.
\end{cases}
\end{equation}

\item For $\nu_j(i)=0$, which means $m=1$, $i=0$, $j=1$, set (cf. \cite[(6.7)]{MR528965})
\begin{equation}
h_{\mu_{1}(0)}(r_1,r_2)=
\begin{cases}
-\log r_2, & r_1\leq r_2,
\\
-\log r_1, & r_2\leq r_1.
\end{cases}
\end{equation}
\end{enumerate}
Then the Green operator  $\mathcal{G}_i^{F}(r_1,y_1,r_2,y_2)$ for co-closed $i$-forms, compactly supported on $C(N)$ and of the type
$g(r)\phi(y)$ with $\phi(y)\in\ker\widetilde{\delta^F}\cap\Omega^{i}(N,F)$, is given by (cf. \cite[(6.10)]{MR528965})
\begin{equation}\label{a5}
\mathcal{G}_i^{F}(r_1,y_1,r_2,y_2)=\sum_{j=1}^{\infty}h_{\mu_j(i)}(r_1,r_2)\phi^{i}_j(y_1)\otimes\phi^{i}_j(y_2).
\end{equation}

In order to obtain the Green operator for forms on $C_{(0,1]}(N)$ which satisfy
either absolute boundary condition at $\{1\}\times N$ or relative boundary condition at $\{1\}\times N$, one must modify $\mathcal{G}_i^{F}(r_1,y_1,r_2,y_2)$ as in \cite[p. 290-291]{MR528965}.

We first consider the case of absolute boundary condition.
\begin{enumerate}[{\rm (I)}]
\item For $\nu_j(i)>0$,
\begin{enumerate}[{\rm (i)}]
\item If $\mu_j(i)>0$, set (compare with \cite[(6.11)]{MR528965})
\begin{equation}
\leftidx{_{\rm a}}h_{\mu_j(i)}(r_1,r_2)=h_{\mu_j(i)}(r_1,r_2)-\frac{a_j^-(i)}{2\nu_j(i)\,a_j^+(i)}(r_1r_2)^{a_j^+(i)}.
\end{equation}

\item If $\mu_j(i)=0$, then $i=0$, $j=1$, set (compare with \cite[(6.17)]{MR528965})
\begin{equation}\label{f2}
\leftidx{_{\rm a}}h_{\mu_1(0)}(r_1,r_2)=h_{\mu_1(0)}(r_1,r_2)+\tfrac{1}{2}(r_1^2+r^2_2)+\tfrac{(1+m)^2}{(1-m)(3+m)}\,.
\end{equation}
\end{enumerate}

\item For $\nu_j(i)=0$, set (compare with \cite[(6.20)]{MR528965})
\begin{equation}\label{f3}
\begin{split}
\leftidx{_{\rm a}}h_{\mu_1(0)}(r_1,r_2)=h_{\mu_1(0)}(r_1,r_2)+\tfrac{1}{2}(r_1^2+r^2_2)-\tfrac{3}{4}\,.
\end{split}
\end{equation}
\end{enumerate}
Then the Green operator  $\mathcal{G}_{{\rm a},i}^{F}(r_1,y_1,r_2,y_2)$ for co-closed $i$-forms on $C_{(0,1]}(N)$
with absolute boundary condition and of the type $g(r)\phi(y)$ is given by
\begin{equation}\label{f1}
\mathcal{G}_{{\rm a},i}^{F}(r_1,y_1,r_2,y_2)=\sum_{j=1}^\infty\leftidx{_{\rm a}}h_{\mu_j(i)}(r_1,r_2)\phi^{i}_j(y_1)\otimes\phi^{i}_j(y_2).
\end{equation}

We now consider the case of relative boundary condition.
\begin{enumerate}[{\rm (I)}]
\item For $\nu_j(i)>0$, set
(compare with \cite[(6.12)]{MR528965})
\begin{equation}
\leftidx{_{\rm r}}h_{\mu_j(i)}(r_1,r_2)=h_{\mu_j(i)}(r_1,r_2)-\frac{1}{2\nu_j(i)}(r_1r_2)^{a_j^+(i)}.
\end{equation}
\item For $\nu_j(i)=0$, set (compare with \cite[(6.22)]{MR528965})
\begin{equation}
\leftidx{_{\rm r}}h_{\mu_1(0)}(r_1,r_2)=h_{\mu_1(0)}(r_1,r_2).
\end{equation}
\end{enumerate}
Then the Green operator  $\mathcal{G}_{{\rm r},i}^{F}(r_1,y_1,r_2,y_2)$ for co-closed $i$-forms on $C_{(0,1]}(N)$
with relative boundary condition and of the type $g(r)\phi(y)$ is given by
\begin{equation}
\mathcal{G}_{{\rm r},i}^{F}(r_1,y_1,r_2,y_2)=\sum_{j=1}^\infty\leftidx{_{\rm r}}h_{\mu_j(i)}(r_1,r_2)\phi^{i}_j(y_1)\otimes\phi^{i}_j(y_2).
\end{equation}

Now let $\dist(x_1,x_2)$ denote the distance from $x_1$ to $x_2$ on $C_{(0,1]}(N)$.
Let ${G}_{{\rm a},i}^{F}(r_1,y_1,r_2,y_2)$ denote the full Green operator on $\Omega^i(C_{(0,1]}(N),F)$ with
absolute boundary condition.
\begin{thm}\label{f6}{\rm (Compare with \cite[Theorems 6.24, 6.43]{MR528965})} There exists a constant $C>0$
such that for $x_1=(r_1,y_1), x_2=(r_2,y_2)\in C_{(0,1]}(N)$,
\begin{align}\label{a4}
\|{G}_{{\rm a},i}^{F}(x_1,x_2)\|\leq C\cdot
\begin{cases}
1+|\log \dist(x_1,x_2)|, &\ m=1,
\\
\dist^{1-m}(x_1,x_2), &\ m\geq 2,
\end{cases}
\end{align}
where by $\|\cdot\|$ we mean the pointwise norm on $C_{(0,1]}(N)$. We will use this notation in the remaining part without further notice.
\end{thm}
\noindent{\bf Proof}\hspace{8pt}Without loss of generality, we assume $r_1\leq r_2$.

We rewrite the co-exact part of \eqref{a5} as follows (cf. \cite[(6.25)]{MR528965}),
\begin{equation}\label{a6}
\begin{split}
\underline{\mathcal{G}}_i^{F}(x_1,x_2):
&=\frac{1}{2}(r_1r_2)^{\alpha(i)}\sum_{\mu_j(i)>0}\frac{(r_1/r_2)^{\nu_j(i)}}{\nu_j(i)}\phi^{i}_j(y_1)\otimes\phi^{i}_j(y_2)
\\
&=\frac{1}{2}(r_1r_2)^{\alpha(i)}P_{{\rm \widetilde{ce}},i}
\frac{e^{\log(r_1/r_2)\cdot(\alpha^2(i)+\widetilde{\Delta})^{1/2}}}{(\alpha^2(i)+\widetilde{\Delta})^{1/2}},
\end{split}
\end{equation}
where by $P_{{\rm \widetilde{ce}},i}$ we mean the orthogonal projection on $\widetilde{\delta^F}\Omega^{i+1}(N,F)$.

Observe that the pointwise norm of \eqref{a6},
viewed as a kernel on $C_{(0,1]}(N)$, is equal to the pointwise norm of \eqref{a6} viewed as a kernel on $N$, multiplied by $(r_1r_2)^{-i}$.

For the $\tfrac{r_1}{r_2}\geq \tfrac{1}{2}$ case, using \cite[Theorem 6.23(1)]{MR528965}, one sees
there exists a constant $C>0$ such that
\begin{align}
&\|\underline{\mathcal{G}}_i^{F}(x_1,x_2)\|\leq C(r_1r_2)^{\frac{1-m}{2}}
\notag
\\
&\hspace{50pt}\cdot
\begin{cases}\label{b1}
1+\big|\log\big (\log^2(r_1/r_2)+\dist^2(y_1,y_2)\big)\big|, &\ m=1,
\\
1+\big(\log^2(r_1/r_2)+\dist^{2}(y_1,y_2)\big)^{\frac{1-m}{2}}, &\ m\geq 2,
\end{cases}
\end{align}
where $\dist(y_1,y_2)$ denote the distance from $y_1$ to $y_2$ on $(N,g^{TN})$.

As in \cite[(6.28)-(6.32)]{MR528965}, one has
\begin{align}
(r_1r_2)^{\frac{1-m}{2}}\leq 2^{3(1-m)/2}\dist^{1-m}(x_1,x_2),
\\
r_1\,r_2\big(\log^2(r_1/r_2)+\dist^2(y_1,y_2)\big)\geq \frac{r_1}{r_2}\dist^2(x_1,x_2),
\\
\label{b2}
\log^2(r_1/r_2)+\dist^2(y_1,y_2)\geq \frac{1}{r^2_2}\dist^2(x_1,x_2).
\end{align}
From \eqref{b1}-\eqref{b2}, one deduces that when $\tfrac{r_1}{r_2}\geq \tfrac{1}{2}$,
\begin{align}
\|\underline{\mathcal{G}}_i^{F}(x_1,x_2)\|\leq C\cdot
\begin{cases}
1+|\log \dist(x_1,x_2)|, &\ m=1,
\\
\dist^{1-m}(x_1,x_2), &\ m\geq 2.
\end{cases}
\end{align}

For the $\tfrac{r_1}{r_2}\leq \tfrac{1}{2}$ case, using Sobolev inequality, one knows the series \eqref{a6},
viewed as a kernel on $N$, converges uniformly. Therefore,
\begin{equation}
\|\underline{\mathcal{G}}_i^{F}(x_1,x_2)\|\leq C(r_1r_2)^{\frac{1-m}{2}}(r_1/r_2)^{\widehat{\nu}(i)},
\ \text{with}\ \widehat{\nu}(i)=\min_{\mu_j(i)>0}\{\nu_j(i)\}.
\end{equation}
From the fact that (cf. \cite[(6.34)]{MR528965})
\begin{equation}
\min_{\mu_j(i)>0}\{\mu_j(i)\}\geq (m-i)(i+1),
\end{equation}
one knows $\widehat{\nu}(i)\geq \tfrac{m+1}{2}$. Thus,
\begin{equation}
\|\underline{\mathcal{G}}_i^{F}(x_1,x_2)\|\leq C r_1r_2^{-m}\leq \tfrac{1}{2}Cr_2^{1-m}\leq 2^{m-2}C \dist^{1-m}(x_1,x_2).
\end{equation}

Now from \eqref{f1}, one sees
\begin{equation}
\begin{split}\label{f4}
\mathcal{G}_{{\rm a},i}^{F}(x_1,x_2)
&=\underline{\mathcal{G}}_{i}^{F}(x_1,x_2)
+\leftidx{_{\rm a}}h_{\mu_1(0)}(r_1,r_2)\phi^{0}_1(y_1)\otimes\phi^{0}_1(y_2)
\\
&\hspace{20pt}-\sum_{\mu_j(i)>0}\frac{a_j^-(i)}{2\nu_j(i)\,a_j^+(i)}(r_1r_2)^{a_j^+(i)}\phi^{i}_j(y_1)\otimes\phi^{i}_j(y_2)\,.
\end{split}
\end{equation}
The last term in \eqref{f4} is the same as the term in \cite[(6.38)]{MR528965}, so
it satisfies the required estimates. On the other hand, it is straightforward to show
the second term in \eqref{f4} implies the required estimates. Thus, one has
\begin{align}\label{f5}
\|\mathcal{G}_{{\rm a},i}^{F}(x_1,x_2)\|\leq C\cdot
\begin{cases}
1+|\log \dist(x_1,x_2)|, &\ m=1,
\\
\dist^{1-m}(x_1,x_2), &\ m\geq 2.
\end{cases}
\end{align}
A similar argument shows that $\mathcal{G}_{{\rm r},i}^{F^*}$ satisfies the estimates
in \eqref{f5}. Moreover, with the help of \cite[Theorem 6.23(2)]{MR528965}, in the same way one can
show that $d^{F}_{x_1}\mathcal{G}_{{\rm a}}^{F}(x_1,x_2)$, $d^{F}_{x_2}\mathcal{G}_{{\rm a}}^{F}(x_1,x_2)$,
$d^{F^*}_{x_1}\mathcal{G}_{{\rm r}}^{F^*}(x_1,x_2)$, $d^{F^*}_{x_2}\mathcal{G}_{{\rm r}}^{F^*}(x_1,x_2)$ are bounded
by $\dist^{-m}(x_1,x_2)$. Since
\begin{equation}
\begin{split}
{G}_{{\rm a},i}^{F}&=\mathcal{G}_{{\rm a},i}^{F}
+d^{F}_{1}\mathcal{G}_{{\rm a},i-1}^{F}\circ d^{F}_{2}\mathcal{G}_{{\rm a},i-1}^{F}
+(*_1^F)^{-1}\mathcal{G}_{{\rm r},m+1-i}^{F^*}*_2^{F}
\\
&\hspace{30pt}+(*_1^F)^{-1}d^{F^*}_{1}\mathcal{G}_{{\rm r},m-i}^{F^*}\circ d^{F^*}_{2}\mathcal{G}_{{\rm r},m-i}^{F^*}*_2^{F}\,,
\end{split}
\end{equation}
applying \cite[Lemma 5.6]{MR528965}, we complete the proof .

\

The main result of this subsection is the Sobolev estimates on $C_{(0,1]}(N)$ as follows.

\begin{thm}\label{d1}
For $n>\frac{m+1}{4}$, there exists a constant $C(n)>0$ with the following property.
Let $\beta\in \Omega^i(C_{(0,1]}(N),F)$ such that $\Delta^j \beta$ satisfies absolute
boundary condition for $0\leq j\leq n$. Then,
\begin{enumerate}[{\rm (i)}]
\item If $i\neq 0$,
\begin{equation}\label{f7}
\|\beta\|\leq C(n) \|\Delta^n \beta\|_{L^2(C_{(0,1]}(N))}\,.
\end{equation}

\item If $i=0$,
\begin{equation}
\|\beta\|\leq C(n) \Big(\|\beta\|_{L^2(C_{(0,1]}(N))}+\|\Delta^n \beta\|_{L^2(C_{(0,1]}(N))}\Big)\,.
\end{equation}
\end{enumerate}
Here by $\|\cdot\|_{L^2(C_{(0,1]}(N_p))}$ we mean the $L^2$ norm on the cone $C_{(0,1]}(N_p)$.
\end{thm}
\noindent{\bf Proof}\hspace{8pt}If $i\neq 0$,
\begin{equation}
\beta(x_1)=({G}_{{\rm a},i}^{F})^n(x_1,x_2)\Delta^n \beta(x_2).
\end{equation}
By \cite[Lemma 5.6]{MR528965} and Theorem \ref{f6}, the norm of $({G}_{{\rm a},i}^{F})^n$ as a function
of $x_2$ is finite provided $n>\frac{m+1}{4}$. Thus, applying the Schwartz inequality, one gets \eqref{f7}.

If $i\neq 0$, the argument is the same expect for the fact that the harmonic function must
be split off and treated separately. We complete the proof.

\

\subsection{Proof of the main result}
Let $(X,g^{TX})$ be a Riemannian manifold of dimension $m+1$ with isolated conical singularities $\Sigma$.
We assume $m$ is odd, and $(X,g^{TX})$ carries an orbifold structure.
As pointed out in \cite[p. 2]{MR2359514}, any orbifold has an atlas consisting of linear charts. Thus,
near each singularity $p\in \Sigma$, we may take the chart of the form
\begin{equation*}
(G_p\subset O(m+1), B^X(p,\varepsilon_p), B(0,\varepsilon_p)),
\end{equation*}
where $B^X(p,\varepsilon_p)$ is the open ball in $X$ of center $p$ and radius $\varepsilon_p>0$,
and $B(0,\varepsilon_p)$ is the open ball in $\mathbb{R}^{m+1}$ of center $0$ and radius $\varepsilon_p$. 
Furthermore, the finite group $G_p$ acts on $B(0,\varepsilon_p)\backslash\{0\}$ freely. 

From our assumption, one knows each $N_p$ ($p\in \Sigma$) in \eqref{t3} is
given by $S^{m}/G_p$, the quotient space induced from a free $G_p$ action on $S^{m}$.

Let $F$ be a flat real orbifold vector bundle over $X$ with the flat connection $\nabla^F$ and an Euclidean metric $g^F$.
Using the identification \eqref{a1} near $p\in \Sigma$, one knows
there exists a $G_p$-equivariant flat vector bundle
$\widetilde{\pi}_p: (\widetilde{F}_p,\nabla^{\widetilde{F}_p},g^{\widetilde{F}_p})\rightarrow B(0,\varepsilon_p)$
such that ${F}\big|_{B^X(p,\varepsilon_p)}\simeq\widetilde{F}_p/G_p$. Since $B(0,\varepsilon_p)$ is contractible,
$\widetilde{F}_p$ is a trivial flat bundle.

We make the assumption that $F$ is an honest vector bundle on $X$. Then the $G_p$-action is trivial along the fiber under our trivialization.
It is this crucial fact that we use in the calculation of the cohomology of $X$ with $F$-coeffcients and
that enables us to apply the result obtained in the last subsection to get the estimate for the heat kernel $K_{\rm c}$ on $X$.

The following theorem compares the $L^2$-cohomology with the singular cohomology.
\begin{thm}\label{ker}
We have the isomorphisms as follows,
\begin{equation}\label{t5}
H_{(2)}^i(X;F)\simeq H^{i}(X;F),\quad 0\leq i\leq m+1.
\end{equation}
Moreover, we have the explicit identifications,
\begin{align}\label{t7}
\ker \Delta_{{\rm c},i} \simeq \ker \Delta_{{\rm o},i}\,,\quad0\leq i\leq m+1.
\end{align}
\end{thm}
\noindent{\bf Proof}\hspace{8pt} Since $X$ has a good cover,
using Mayer Vietoris sequence, it is sufficient to
show for $p\in \Sigma$ and $u>0$ small enough,
\begin{equation}\label{y1}
H_{(2)}^i(C^*_{(0,u]}(N_p);F)\simeq H^{i}(C^*_{(0,u]}(N_p);F),\quad 0\leq i\leq m+1.
\end{equation}

We have the following fact from \cite{MR573430},
\begin{align}\label{y2}
H_{(2)}^i(C^*_{(0,u]}(N_p);F)\simeq
\begin{cases}
H^i(N_p;F),    &i\leq \frac{m}{2},
\\0, &i\geq \frac{m+1}{2}.
\end{cases}
\end{align}
On the other hand,
\begin{align}\label{y3}
H^i(N_p;F)\simeq \big(H^{i}(S^{m};\widetilde{F}_p)\big)^{G_p}\simeq
\begin{cases}
F|_{\{p\}},    &i=0, \text{or}\ i=m,
\\0, &i\neq 0,\text{and}\ i\neq m.
\end{cases}
\end{align}
Since $C^*_{(0,u]}(N_p)$ is contractible, \eqref{y1} follows from \eqref{y2} and \eqref{y3}.

Observe that for any element $s\in \ker\Delta_{{\rm o},i}$,
$s$ and $d^Fs$ are both $L^2$ integrable, hence $s$ determines an element in $\ker\Delta_{{\rm c},i}$.
From \eqref{c1}, \eqref{c2} and \eqref{t5}, one sees the injection from $\ker\Delta_{{\rm o},i}$ to $\ker\Delta_{{\rm c},i}$
gives the isomorphism in \eqref{t7}.

\begin{rem}
One can verifies {\rm \eqref{t7}} directly. In fact,
using the explicit expression of the harmonic forms on the cone in {\rm \cite{MR528965,MR730920}}, one finds the harmonic element in
$\ker\Delta_{\rm c}$ has removable singularities, so it lies in $\ker\Delta_{\rm o}$ by elliptic regularity.
\end{rem}

\begin{thm}\label{es1}
Given $T>0$, $0<u_0\leq 1$, $n\in \mathbb{N}$, and $n>\frac{m+1}{4}$, there exists a constant $C(T,u_0,n)>0$ such that
for $0<u<\frac{1}{4}u_0$, $x\in X_u$, $y\in X_{u_0}$, and $0<t\leq T$,
\begin{equation}
\|{K}_{\rm c}(t,x,y)\|\leq C(T,u_0,n)\,t^n.
\end{equation}
The same estimates hold for $d_x^F{K}_{\rm c}(t,x,y)$ and $\delta_y^F {K}_{\rm c}(t,x,y)$.
\end{thm}
\noindent{\bf Proof}\hspace{8pt}As explained in \cite[p. 287]{MR528965}, it suffices to
prove the estimate holds for ${K}_{\rm c}(t,x,y)$.
We first establish the estimates for the
heat kernel $\widehat{K}_{p}(t,x,y)$ on the bounded cone ${C_{(0,1]}}(N_p)$.

Choose a parametrix $P_n(t,x,y)$ of order $n$ \cite[p. 272]{MR528965} such that
\begin{equation}
P_n(0,x,y)=\delta_y, \ \text{and}\ P_n(t,x,y)=0,\ \text{for}\ x\in X_u, y\in X_{u_0}.
\end{equation}
Then there exists a constant $C'(T,u_0,n)>0$ such that for $j\in\mathbb{N}$, $t\in (0,T]$
\begin{equation}\label{d2}
\big\|\Delta^j\big(\widehat{K}_{p}(t,x,y)-P_n(t,x,y)\big)\big\|_{L^2(C_{(0,1]}(N_p))}\leq C'(T,u_0,n)\,t^n,
\end{equation}
where by $\|\cdot\|_{L^2(C_{(0,1]}(N_p))}$ we mean taking the pointwise norm with respect to the $x$ variable
and the $L^2$ norm with respect to the $y$ variable with $x$, $y$ varying on the cone $C_{(0,1]}(N_p)$.

By \eqref{d2} and the Sobolev inequality on the bounded cone, Theorem \ref{d1},
one sees for $n>\frac{m+1}{4}$, there exists a constant $C(T,u_0,n)>0$ such that the pointwise norm satisfies
\begin{equation}\label{d3}
\big\|\widehat{K}_{p}(t,x,y)-P_n(t,x,y)\big\|\leq C(T,u_0,n)\,t^n.
\end{equation}
In particular, from \eqref{d3} one gets for $x\in X_u, y\in X_{u_0}$,
\begin{equation}
\big\|\widehat{K}_{p}(t,x,y)\big\|\leq C(T,u_0,n)\,t^n.
\end{equation}
As explained in \cite[p. 286]{MR528965}, the same estimates hold
for $\frac{\partial^\alpha}{\partial y^\alpha}\widehat{K}_p(t,x,y)$.

Let $\phi(r_2):[0,1]\rightarrow[0,1]$ be a smooth
function such that
\begin{equation}
\phi(r_2)=
\begin{cases}
1, &0\leq r_2\leq \frac{1}{3}u_0,
\\
0, &\frac{2}{3}u_0\leq r_2\leq 1.
\end{cases}
\end{equation}
Then as in \cite[(7.7)]{MR528965}, one has  for $j\in\mathbb{N}$,
\begin{align}
&\|\Delta_y^j{K}_{\rm c}(t,x,y)\|_{L^2(X_{u_0})}
\notag\\
&\leq\Big\|\Delta_y^j\Big({K}_{\rm c}(t,x,y)-\sum_{p\in \Sigma}\phi(r_2)\widehat{K}_p(t,x,y)\Big)\Big\|_{L^2(X\backslash\Sigma)}
\notag\\
&\leq\sum_{p\in \Sigma}\int_{0}^t\Big\|\Delta_y^j\Big(\partial_s+\Delta_y\Big)\Big(\phi(r_2)\widehat{K}_p(s,x,y)\Big)\Big\|_{L^2(C_{(0,1]}(N_p))}\ ,
\end{align}
where the norm $\|\cdot\|_{L^2}$ is the pointwise norm with respect to the $x$ variable
and the $L^2$ norm with respect to the $y$ variable as in \eqref{d2}.

In view of the estimates already
established for $\frac{\partial^\alpha}{\partial y^\alpha}\widehat{K}_p(t,x,y)$, the estimates
for ${K}_{\rm c}(t,x,y)$, $d_x^F{K}_{\rm c}(t,x,y)$ and $\delta_y^F {K}_{\rm c}(t,x,y)$ follow
from the standard Sobolev inequality (cf. \cite[(5.7)]{MR528965}) applied with respect to the $y$ variable.

\

We are now in a position to prove our main theorem which compares the two heat kernels outside the singular points.
\begin{thm}\label{g3}
For $t>0$, $x,y\in X\backslash\Sigma$, one has
\begin{equation}\label{t1}
{K}_{{\rm o},i}(t,x,y)={K}_{{\rm c},i}(t,x,y),\quad 0\leq i\leq m+1.
\end{equation}
\end{thm}
\noindent{\bf Proof}\hspace{8pt}Fix any $T>0$, and we will show \eqref{t1} holds for any $t\in (0,T]$.

For any $u_0>0$ small enough fixed temporarily and $x,y\in X_{u_0}$, applying the Duhamel principle \cite[(3.9)]{MR528965} on $X_u$ with
$0<u<\frac{1}{4}u_0$, one gets
\begin{align}
K_{{\rm c},i}(t,x,y)-K_{{\rm o},i}(t,x,y)
&=\int_0^t\int_{\partial X_u}\big\langle K_{{\rm o},i}(t-s,x,z)\wedge*^F d^F K_{{\rm c},i}(s,z,y)\big\rangle
\notag
\\
&\hspace{-50pt}-\int_0^t\int_{\partial X_u}\big\langle*^Fd^F K_{{\rm o},i}(t-s,x,z)\wedge K_{{\rm c},i}(s,z,y)\big\rangle
\notag
\\
&\hspace{-60pt}+\int_0^t\int_{\partial X_u}\big\langle\delta^F K_{{\rm o},i}(t-s,x,z)\wedge *^F K_{{\rm c},i}(s,z,y)\big\rangle
\notag\\
&\hspace{-70pt}-\int_0^t\int_{\partial X_u}\big\langle *^F K_{{\rm o},i}(t-s,x,z)\wedge\delta^F K_{{\rm c},i}(s,z,y)\big\rangle,
\end{align}
where on the right hand side all operations are applied to the variable $z$.

Using the uniform estimates in Theorems \ref{es2}, \ref{es1}, one sees  as $u\rightarrow 0$,
\begin{align}\label{t2}
K_{{\rm o},i}(t,x,y)=K_{{\rm c},i}(t,x,y)\text{ holds for }x, y \in X_{u_0}.
\end{align}
Since $u_0>0$ can be taken arbitrarily small, \eqref{t2} holds on $X\backslash\Sigma$.

\

From Theorems \ref{ker}, \ref{g3}, one gets

\begin{cor}The two torsions are equal, i.e.,
\begin{align}
T_{\rm c}(X,g^{TX},g^F)= T_{\rm o}(X,g^{TX},g^F).
\end{align}
\end{cor}


\end{document}